\documentclass{elsart}
\def\Z/{{\mathbb Z}}
\def\Zn/{{\mathbb Z}_n}
\def\beh/{behaving}
\def\zf/{zero-free}
\def\zfs/{zero-free sequence}
\def\pin/{positive integer}
\def\pins/{positive integer sequence}
\def\lpr/{least positive representatives}
\def\v#1{\overline{#1}}
\def\length/{length greater than~$n/2$}
\def\lengthell/{length~$\ell>n/2$}
\def\threshold/{(2n{-}2)/3}
\def\bound/{\ell-\left\lfloor\frac{n{-}1}{3}\right\rfloor}
\def\boundnew/{\ell-\lfloor(n{-}1)/3\rfloor}
\newtheorem{quest}[thm]{Question}
\usepackage{amssymb,amsmath}

\begin{document}

\begin{frontmatter}
\title{Long \zfs/s in finite cyclic groups}
\author[]{Svetoslav Savchev}
\kern-.35em\footnote{No current affiliation.},
\author[Emory]{Fang Chen}
\address[Emory]{Oxford College of Emory University, Oxford, GA 30054, USA}

\begin{abstract}
A sequence in an additively written abelian group is called \zf/ if
each of its nonempty subsequences has sum different from the zero
element of the group. The article determines the structure of the
\zfs/s with
lengths greater than~$n/2$ in the additive group~$\Zn/$ of integers
modulo~$n$. The main result states that for each
\zfs/~$(a_i)_{i=1}^\ell$ of \lengthell/ in~$\Zn/$ there is an
integer~$g$ coprime to~$n$ such that if $\v{ga_i}$ denotes the least
\pin/ in the congruence class~$ga_i$ (modulo~$n$), then
$\Sigma_{i=1}^\ell\v{ga_i}<n$. The answers to a number of frequently
asked zero-sum questions for cyclic groups follow as immediate
consequences. Among other applications, best possible lower bounds
are established for the maximum multiplicity of a term in a \zfs/
with \length/, as well as for the maximum multiplicity of a
generator. The approach is combinatorial and does not appeal to
previously known nontrivial facts.

\end{abstract}

\end{frontmatter}

\section{Introduction}
\label{Introduction}

Among $n$~arbitrary integers one can choose several whose sum is
divisible by~$n$. In other words, each sequence of length~$n$ in the
cyclic group of order~$n$ has a nonempty subsequence with sum~zero.
This article describes all sequences of \length/ in the same group
that fail the above property.

Here and henceforth, $n$ is a fixed integer greater
than~1, and the cyclic group of order~$n$ is identified with the
additive group~$\Zn/=\Z//n\Z/$ of integers modulo~$n$. A sequence
in~$\Zn/$ is called a {\em zero sequence} or a {\em zero sum} if the
sum of its terms is the zero element of~$\Zn/$. A sequence is {\em
\zf/} if it does not contain nonempty zero subsequences.

We study the general structure of the \zfs/s in~$\Zn/$ whose lengths
are between~$n/2$ and~$n$.
Few nontrivial related results
are known to us, of which we mention only one. A work of
Gao~\cite{Gao} characterizes the \zfs/s of length roughly greater
than $2n/3$. On the other hand, structural information about shorter
\zfs/s
naturally translates into knowledge about
problems
of
significant
interest. Several examples to this effect are
included below. The main result provides complete answers to a
number of repeatedly explored zero-sum questions.

Our objects of study
can be characterized in very
simple terms. To be more specific, let us recall several standard
notions.

If $g$ is an integer coprime to~$n$, multiplication by~$g$ preserves
the zero sums in~$\Zn/$ and does not introduce new ones. Hence a
sequence $\alpha=(a_1,\dots,a_k)$ is \zf/ if and only if the
sequence $g\alpha=(ga_1,\dots,ga_k)$ is \zf/, which motivates the
following definition.

For sequences $\alpha$ and $\beta$ in~$\Zn/$, we say that $\alpha$
is {\em equivalent} to~$\beta$ and write $\alpha\cong \beta$ if
$\beta$ can be obtained from~$\alpha$ through multiplication by an
integer coprime to~$n$ and rearrangement of terms. Clearly $\cong$
is an equivalence relation.

If $\alpha=(a_1,\dots,a_k)$ is a sequence in~$\Zn/$, let $\v{a_i}$
be the unique integer in the set~$\{1,2,\dots,n\}$ which belongs to
the congruence class~$a_i$ modulo~$n$, $i=1,\dots,k$. The
number~$\v{a_i}$ is called the {\em least positive representative}
of~$a_i$. Consequently, the sum $L(\alpha)=\sum_{i=1}^k\v{a_i}$ will
be called the {\em sum of the \lpr/} of~$\alpha$.

Now the key result in the article, Theorem~\ref{major}, can be
stated as follows:
\begin{quote}\em
Each \zfs/ of \length/ in~$\Zn/$ is equivalent to a sequence whose
sum of the \lpr/ is less than~$n$.
\end{quote}
This statement reduces certain zero-sum problems in cyclic groups to
the study of easy-to-describe \pins/s. Thus all proofs in
Sections~\ref{Index}--\ref{NewFunction} are merely short elementary
exercises.

The approach of the article is combinatorial and does not follow a
line of thought known to us from previous work. The exposition is
self-contained in the sense that it does not rely on any nontrivial
general fact. Sections~\ref{Preliminaries}
and~\ref{BehavingSequences} are preparatory. The main result is
proven in~Section~\ref{TheMainResult}.

For a sequence~$\alpha$ in~$\Zn/$, the number $Index(\alpha)$ is
defined as the minimum of $L(g\alpha)$ over all~$g$ coprime to~$n$.
Section~\ref{Index} contains the answer, for all~$n$, to the
question about the minimum~$\ell(\Zn/)$ such that each minimal zero
sequence of length at least~$\ell(\Zn/)$ in~$\Zn/$ has index~$n$.

Issues of  considerable interest among the zero-sum problems are the
maximum multiplicity of a term in a \zfs/, and of a generator in
particular. Sections~\ref{MaxMultTerm} and~\ref{MaxMultGenerator}
provide exhaustive answers for \zfs/s of all lengths~$\ell>n/2$
in~$\Zn/$. Best possible lower bounds are established in both cases,
which improves on earlier work of Bovey, Erd\H{o}s and Niven
\cite{BoveyErdosNiven}, Gao and Geroldinger \cite{GaoGeroldinger},
Geroldinger and Hamidoune \cite{GeroldingerHamidoune}.

In~Section~\ref{NewFunction} we introduce a function closely related
to the \zfs/s in cyclic groups. This is an analogue of a function
defined by Bialostocki and Lotspeich \cite{BialostockiLotspeich} in
relation to the theorem of Erd\H{o}s, Ginzburg and Ziv
\cite{ErdosGinzburgZiv}. Theorem~\ref{major} enables us to determine
the values of the newly defined function in a certain range. An
explicit description of the \zfs/s with a given \lengthell/
in~$\Zn/$ is included in Section~\ref{ConcludingRemarks}.

\section{Preliminaries}
\label{Preliminaries}

Several elementary facts about sequences in general abelian groups
are
considered below. We precede them by remarks on terminology and
notation. The {\em sumset} of a sequence in an abelian group~$G$ is
the set of all $g\in G$ representable as a nonempty subsequence sum.
The cyclic subgroup of~$G$ generated by an element~$g\in G$ is
denoted by~$\langle g\rangle$; the order of $g$ in~$G$ is denoted
by~$\text{ord(g)}$.

\begin{prop}\label{increase}
For a \zfs/ $(a_1,\dots,a_k)$ in an abelian group, let $\Sigma_i$ be
the sumset of the subsequence $(a_1,\dots,a_i)$, $i=1,\dots,k$. Then
$\Sigma_{i-1}$ is a proper subset of $\Sigma_i$ for each
$i=2,\dots,k$. Moreover, the subsequence sum $a_1+\cdots+a_i$
belongs to $\Sigma_i$ but not to $\Sigma_{i-1}$. In particular,
$a_1+\cdots+a_k$ belongs to~$\Sigma_k$ but not to any~$\Sigma_i$
with~$i<k$.
\end{prop}

\begin{pf}
Since $\Sigma_{i-1}\subseteq \Sigma_i$ and $a_1+\cdots+a_i\in
\Sigma_i$, it suffices to prove that $a_1+\cdots+a_i\not\in
\Sigma_{i-1}$, $i=2,\dots,k$. Suppose
that $a_1+\cdots+a_i\in \Sigma_{i-1}$ for some $i=2,\dots,k$. Then
$a_1+\cdots+a_i=\sum_{j\in J}a_j$ for a nonempty subset $J$ of
$\{1,\dots,i-1\}$. Each term on the right-hand side is present on
the left-hand side, and $a_i$ is to be found only on the left. So
canceling yields a nonempty zero sum in $(a_1,\dots,a_k)$, which
contradicts the assumption that it is \zf/.\qed
\end{pf}

Proposition~\ref{increase} states that, for a \zfs/
$\alpha=(a_1,\dots,a_k)$, the sumset of the subsequence
$(a_1,\dots,a_{i{-}1})$  strictly increases upon appending the next
term~$a_i$, $i=2,\dots,k$. If the increase of the sumset size is
exactly~1, we say that $a_i$ is a {\em $1$-term} for~$\alpha$.
Naturally, the property of being a 1-term is not necessarily
preserved upon rearrangement of terms.

The next statement contains observations on 1-terms. Parts~a) and~b)
seem to be folklore and can be found for instance
in~\cite{SmithFreeze}.

\begin{prop}\label{arabian}
Let $\alpha=(a_1,\dots,a_k)$ be a nonempty \zfs/ with sumset
$\Sigma$ in an abelian group~$G$. Suppose that, for some $b\in G$,
the extended sequence $\alpha\cup\{b\}=(a_1,\dots,a_k,b)$  is
\zf/ and $b$ is a $1$-term for~$\alpha\cup\{b\}$. Then:
\begin{enumerate}
\item[a)] $\Sigma$ is the union of a progression $\{b,2b,\dots,sb\}$, where
$1\le s <\text{\em ord}(b){-}1$, and several (possibly none)
complete proper cosets of the cyclic subgroup generated by~$b$;

\item[b)] the sum of $\alpha$ equals $sb$;

\item[c)] $b$ is the unique element of~$G$ that can be appended to~$\alpha$ as a last term
so that the resulting sequence is \zf/ and ends in a $1$-term.
\end{enumerate}
\end{prop}

\begin{pf}
Parts a) and b) are proven in~\cite{SmithFreeze}. For part~c), let
$c\in G$ be such that the
sequence $\alpha\cup\{c\}=(a_1,\dots,a_k,c)$ is \zf/ and $c$ is a
1-term for~$\alpha\cup\{c\}$. We prove that $c=b$. Because $b$ is a
1-term for~$\alpha\cup\{b\}$, in view of a) we have
$\Sigma=\{b,2b,\dots,sb\}\cup C_1\cup\cdots\cup C_m$, where $1\le s
<\text{ord}(b){-}1$ and $C_1,\dots,C_m$ are complete proper cosets
of the subgroup $\langle b\rangle$ generated by~$b$. The sumset
$\Sigma'$ of $\alpha\cup\{c\}$ contains the progression
$P=\{c,c+b,\dots,c+sb\}$ whose length $s{+}1$ is at least~2. Since
$c$ is a 1-term for~$\alpha\cup\{c\}$, it follows that $P$
intersects $\{b,2b,\dots,sb\}$ or one of
$C_1,\dots,C_m$. By b), $P$ contains the sum $c+sb$ of
$\alpha\cup\{c\}$, which is
an element of $\Sigma'\setminus \Sigma$ in view of
Proposition~\ref{increase}. Hence $P\cap C_i=\emptyset$ for all
$i=1,\dots,m$, or else $c+sb\in\Sigma$.
Thus $P$ intersects $\{b,2b,\dots,sb\}$, and $0\not\in P$
implies $c=xb$ for some integer
$x$ satisfying $1\le x\le s$. Hence the progression
$\{b,2b,\dots,(s+x)b\}$ is contained in $\Sigma'$. Now we see that
the size of $\Sigma$ grows exactly by 1 upon appending $c$
only if $x=1$, i.~e. $c=b$.\qed
\end{pf}

A \zfs/ in a finite abelian group~$G$ is {\em maximal} if it is not
a subsequence of a longer \zfs/ in~$G$. Let $\alpha$ be a \zfs/
in~$G$ whose sumset does not contain at least one nonzero
element~$g$ of~$G$. Then $\alpha\cup\{-g\}$ is a longer \zfs/
containing~$\alpha$. This remark and Proposition~\ref{increase} show
that a \zfs/ in~$G$ is maximal if and only if its sumset
is~$G\setminus\{0\}$. The same remark (with
Proposition~\ref{increase} again) yields a quick justification of
the next statement. We omit the proof.

\begin{prop}\label{maximal}
Each \zfs/ in a finite abelian group can be extended to
a maximal \zfs/.
\end{prop}

\section{Behaving sequences}
\label{BehavingSequences}

A \pins/ with sum $S$ will be called {\em \beh/} if its sumset is
$\{1,2,\dots,S\}$. The ordering of the sequence terms is not
reflected in the definition. However, assuming them in nondecreasing
order enables one to state a convenient equivalent description. Its
sufficiency part is a problem from the 1960 edition of the
celebrated K\"ursch\'ak contest in Hungary, the oldest mathematics
competition for high-school students in the world.

\begin{prop}\label{kurschak}
\label{behaving} A sequence $(s_1,\dots,s_k)$ with positive integer
terms in nondecreasing order $s_1\le\cdots\le s_k$ is \beh/ if and
only if
\begin{equation*}
s_1=1\qquad\text{and}\qquad s_{i+1}\le
1+s_1+\cdots+s_i\quad\text{for all $i=1,\dots,k-1$.}
\end{equation*}
\end{prop}

\begin{pf}
Denote $S=s_1+\cdots+s_k$ and suppose that the sequence is \beh/;
then its sumset is $\Sigma=\{1,2,\dots,S\}$. Since $1\in\Sigma$ and
$s_i\ge 1$ for all $i$, it follows that $s_1=1$. For each
$i=1,\dots,k-1$, let $T_i=1+s_1+\cdots+s_i$. Clearly $T_i\le S$,
hence $T_i\in \Sigma$. Also $T_i>s_1+\cdots +s_i$, so the
subsequence whose sum equals $T_i$ contains a summand~$s_j$ with
index~$j$ greater than~$i$. Therefore $T_i\ge s_j\ge s_{i+1}$, as
desired.

Conversely, let $s_1=1$ and $s_{i+1}\le 1+s_1+\cdots+s_i$,
$i=1,\dots,k-1$. Denoting $S_k=s_1+\cdots+s_k$,
we prove by induction on~$k$ that the sumset of $(s_1,\dots,s_k)$ is
$\{1,2,\dots,S_k\}$. The base $k=1$ is clear. For the inductive
step, let $\Sigma_{k-1}$ and $\Sigma_{k}$ be the sumsets of
$(s_1,\dots,s_{k-1})$ and $(s_1,\dots,s_{k-1},s_{k})$, respectively.
Since $\Sigma_{k-1}=\{1,2,\dots,S_{k-1}\}$ by the induction
hypothesis, it follows that
$\Sigma_{k}=\{1,2,\dots,S_{k-1}\}\cup\{s_{k},s_{k}+1,\dots,s_{k}+S_{k-1}\}$.
In view of the condition $s_{k}\le 1+S_{k-1}$, we obtain
$\Sigma_{k}=\{1,2,\dots,s_{k}+S_{k-1}\}=\{1,2,\dots,S_{k}\}$. The
induction is complete.\qed
\end{pf}

A simple consequence of Proposition~\ref{kurschak} proves essential
for the main proof.

\begin{prop}\label{behaving; length>n/2} Let $k$ be a
positive integer. Each sequence with positive integer terms of
length at least~$k/2$ and sum less than~$k$ is behaving.
\end{prop}

\begin{pf}
Denoting the sequence by $(s_1,\dots,s_\ell)$ and assuming
$s_1\le\cdots\le s_\ell$, we check the sufficient condition of
Proposition~\ref{kurschak}. Given that $\ell\ge k/2$ and
$\Sigma_{i=1}^\ell s_i<k$, it is easy to see that $s_1=1$. Suppose
that $s_{i+1}\ge 2+s_1+\cdots+s_i$ for some $i=1,\dots,\ell-1$. Then
$s_j\ge i+2$ for all~$j=i+1,\dots,\ell$. Therefore
\[
k>\Sigma_{i=1}^\ell s_i\ge i+(\ell-i)(i+2)=2\ell+i(\ell-i-1)\ge
k+i(\ell-i-1)\ge k,
\]
which is a contradiction. The claim follows.\qed
\end{pf}

Now we introduce a key notion. Let $G$ be an abelian group and $g$ a
nonzero element of $G$. A sequence $\alpha$ in $G$ will be called
{\em \beh/ with respect to~$g$} or {\em $g$-\beh/} if it has the
form $\alpha=(s_1g,\dots,s_kg)$, where $(s_1,\dots,s_k)$ is a \beh/
positive integer sequence with sum $S=s_1+\cdots+s_k$ less than the
order of $g$ in $G$.

It follows from the definition that $1\le s_i<\text{ord}(g)$ for
$i=1,\dots,k$. All terms of~$\alpha$ are contained in the cyclic
subgroup~$\langle g\rangle$ generated by~$g$. Moreover, since the
sumset of $(s_1,\dots,s_k)$ is $\{1,2,\dots,S\}$, the sumset of
$\alpha$ is the progression $\{g,2g,\dots,Sg\}$ which is entirely
contained in~$\langle g\rangle$. Finally, $g$ is a term of~$\alpha$
by Proposition~\ref{kurschak} as one of $s_1,\dots,s_k$ equals~1.

\section{The main result}
\label{TheMainResult}

The proof of the main theorem involves certain rearrangements of
terms in \zfs/s. The next lemma states a condition guaranteeing that
such rearrangements are possible.

\begin{lem}\label{rearrangement lemma}
Let $\alpha$ be a \zfs/ of length~$\ell$ greater than~$n/2$
in~$\Zn/$. Suppose that, for some $k\in\{1,\dots,\ell-2\}$, the
first $k+1$~terms of~$\alpha$ form a subsequence with sumset of size
at least~$2k+1$. Then the remaining terms of~$\alpha$ can be
rearranged so that the sequence obtained ends in a $1$-term.
\end{lem}

\begin{pf}
Regardless of how the
last $\ell-k-1$~terms of~$\alpha$ are permuted, at least one of them
will be a $1$-term for the permuted sequence. If not, by
Proposition~\ref{increase} each term after the first $k+1$ increases
the sumset size by at least~2. Hence the total sumset size is at
least $(2k+1)+2(\ell-k-1)=2\ell-1\ge n$ which is impossible for a
\zfs/.

Fix the initial $k+1$ terms of~$\alpha$. Choose a rearrangement of
the last ${\ell-k-1}$ terms
such that the first $1$-term among them occurs as late as possible.
Let this term be~$c$, and let $\alpha'$ be the resulting
rearrangement of~$\alpha$. We are done if $c$ is the last term
of~$\alpha'$. If not, interchange~$c$ with any term~$d$ following it
in~$\alpha'$ to obtain a new rearrangement~$\alpha''$. The same
sequence~$\beta$ precedes~$c$ and~$d$ in~$\alpha'$ and~$\alpha''$,
respectively, and $\beta$ contains no $1$-terms after the
initial~$k+1$ terms. On the other hand, by the extremal choice
of~$\alpha'$, a $1$-term must occur among the last $\ell-k-1$~terms
of~$\alpha''$ at the position of~$d$ in the latest. Therefore $d$ is
a $1$-term for~$\alpha''$. Thus if either of~$c$ and~$d$ is appended
to~$\beta$,
the sequence obtained ends in a $1$-term.
Now
Proposition~\ref{arabian}~c) implies $c=d$. Hence the terms
after~$c$ in~$\alpha'$
are all equal to~$c$, so they are all $1$-terms for~$\alpha'$ by
Proposition~\ref{arabian}~a). In particular, $\alpha'$ ends in a
$1$-term.\qed
\end{pf}

\begin{thm}\label{main}
Each \zfs/ of \length/ in the cyclic group~$\Zn/$ is \beh/ with
respect to one of its terms.
\end{thm}

\begin{pf}
First we prove the theorem for maximal sequences. Let $\alpha$ be a
maximal \zfs/ of \lengthell/ in~$\Zn/$.

For each term~$a$ of~$\alpha$ there exist $a$-\beh/ subsequences
of~$\alpha$, for instance the one-term subsequence $(a)$.
We assign to $a$ one such $a$-\beh/
subsequence~$\alpha_a=(s_1a,\dots,s_ka)$ of maximum length~$k$. Here
$(s_1,\dots,s_k)$ is a \beh/ \pins/ such that $S=s_1+\cdots+s_k$ is
less than the order $\text{ord}(a)$ of $a$ in~$\Zn/$. In particular
$1\le s_i<\text{ord}(a)$, $i=1,\dots,k$. The sumset
of~$(s_1,\dots,s_k)$ is $\{1,2,\dots,S\}$, and the sumset
of~$\alpha_a$ is $\{a,2a,\dots,Sa\}$, a progression contained in the
cyclic subgroup~$\langle a\rangle$ generated by~$a$. Observe that
all occurrences of $a$ in~$\alpha$ are terms of~$\alpha_a$.

We show that there is a term~$g$ whose associated $g$-\beh/
subsequence $\alpha_g$ is the entire~$\alpha$. To this end, choose
an arbitrary term~$a$ of~$\alpha$ and suppose that
$\alpha_a\ne\alpha$. The notation for~$\alpha_a$ from the previous
paragraph is assumed. Let us rearrange~$\alpha$ as follows. Write
the terms of~$\alpha_a$ first and then any term~$b$ of~$\alpha$
which is not in~$\alpha_a$. The subsequence
$\alpha_a\cup\{b\}=(s_1a,\dots,s_ka,b)$ obtained so far has sumset
$P_1\cup P_2$ where $P_1=\{a,2a,\dots,Sa\}$
and~$P_2=\{b,b+a,\dots,b+Sa\}$.

It is not hard to check that $P_1\cap P_2=\emptyset$. This is clear
if $b\not\in\langle a\rangle$ as $P_1$ and~$P_2$ are in different
cosets of~$\langle a\rangle$. Let $b\in \langle a\rangle$, so $b=sa$
with $1\le s<\text{ord}(a)$. Then $P_2=\{sa,(s+1)a,\dots,(s+S)a\}$
and it suffices to prove the inequalities $S+1<s$ and
$s+S<\text{ord}(a)$.

First, $s+S\ge \text{ord}(a)$ implies that $\text{ord}(a)$ occurs
among the consecutive integers $s,s+1,\dots,s+S$.
Hence $P_2$ contains the zero element of~$\Zn/$ which is false.
Next, suppose that $s\le S+1$. Then the integer sequence
$(s_1,\dots,s_k,s)$ has sum $s+S$ and sumset
$\{1,\dots,S,\dots,s+S\}$, so it is \beh/. We also have
$s+S<\text{ord}(a)$, as just shown. But then
$\alpha_a\cup\{b\}=(s_1a,\dots,s_ka,sa)$ is an $a$-\beh/ subsequence
of~$\alpha$ longer than~$\alpha_a$, contradicting the maximum choice
of~$\alpha_a$. Therefore
$P_1$ and~$P_2$ are disjoint also in the case $b\in\langle
a\rangle$.

Now, $P_1\cap P_2=\emptyset$ and $|P_1|=S\ge k$, $|P_2|=S+1\ge k+1$
imply that $|P_1\cup P_2|\ge 2k+1$. It also follows that there are
terms of~$\alpha$ out of $\alpha_a\cup\{b\}$. Otherwise $k+1=\ell$
and because $n-1\ge |P_1\cup P_2|\ge 2k+1$ ($\alpha_a\cup\{b\}$ is
\zf/, hence its sumset has size at most $n-1$), we obtain $n\ge
2\ell$ which is not the case.
Therefore, by Lemma~\ref{rearrangement lemma}, the terms of~$\alpha$
not occurring in~$\alpha_a\cup\{b\}$ can be permuted to obtain a
rearrangement $\alpha'$ which ends in a 1-term~$c$.

Recall now that $\alpha$ is maximal, and hence so is its
rearrangement $\alpha'$. Let $\Sigma$ be the sumset of the sequence
obtained from $\alpha'$ by deleting its last term~$c$. Since $c$ is
a 1-term for~$\alpha'$, $\Sigma$ is missing exactly one nonzero
element of~$\Zn/$. By Proposition~\ref{increase}, the missing
element is the sum~$A\ne 0$ of all terms of~$\alpha$.
On the other hand, $\Sigma$ must be missing the element
$-c$ of~$\Zn/$ ($-c\ne 0$), or else appending~$c$ to
obtain~$\alpha'$ would produce a zero sum. Because the missing
element is unique, we obtain $A=-c$, i.~e. $c=-A$.

We reach the following conclusion. If $\alpha_a\ne\alpha$ for at
least one term~$a$ of~$\alpha$ then the group element~$-A$ is a term
of~$\alpha$. Moreover, if $a$ is any term such that
$\alpha_a\ne\alpha$, the subsequence~$\alpha_a$ does not contain at
least one occurrence of~$-A$.

Apply this conclusion to an arbitrary term~$g$ of~$\alpha$. The
statement is proven if $\alpha_g=\alpha$. If not then $h=-A$ is a
term of~$\alpha$. Consider its associated maximal $h$-\beh/
subsequence~$\alpha_h$. Since $\alpha_h$ contains all occurrences
of~$-A=h$, it follows that $\alpha_h=\alpha$. This completes the
proof in the case where $\alpha$ is maximal.

Suppose now that $\alpha$ is not maximal. By
Proposition~\ref{maximal}, it can be extended to a maximal
\zfs/~$\beta$ in~$\Zn/$, of length~$m>\ell>n/2$. (Clearly $m<n$.) By
the above, there is a term~$a$ of $\beta$ such that $\beta$ is
$a$-\beh/. This is to say, $\beta=(s_1a,\dots,s_{m}a)$  for some
\beh/ positive integer sequence $(s_1,\dots,s_{m})$ with sum less
than $\text{ord}(a)$. Deleting the additionally added terms
from~$\beta$, we infer that $\alpha=(s_{i_1}a,\dots,s_{i_\ell}a)$
for some positive integer sequence $(s_{i_1},\dots,s_{i_\ell})$ of
length~$\ell$ and sum less than $\text{ord}(a)$. Now, since
$\ell>n/2\ge \text{ord}(a)/2$, one can apply
Proposition~\ref{behaving; length>n/2} with $k=\text{ord}(a)$, which
shows that $(s_{i_1},\dots,s_{i_\ell})$ is \beh/. Hence
$\alpha=(s_{i_1}a,\dots,s_{i_\ell}a)$ is $a$-\beh/. Also, $a$ is a
term of~$\alpha$:
as already explained, one of the integers $s_{i_1},\dots,s_{i_\ell}$
equals~$1$ by Proposition~\ref{behaving}. The proof is complete.
\qed
\end{pf}

By Theorem~\ref{main}, each \zfs/ of \lengthell/ in~$\Zn/$ has the
form $\alpha=(s_{1}a,\dots,s_{\ell}a)$, where $a$ is one of its
terms and $(s_1,\dots,s_\ell)$ is a \pins/ with sum less than
$\text{ord}(a)$. In particular $1\le s_i<\text{ord}(a)$ for
$i=1,\dots,\ell$. It is immediate that $\text{ord}(a)=n$. Otherwise
the subgroup $\langle a\rangle$, of order at most $n/2$, would
contain a \zfs/ of \lengthell/ which is impossible. Hence there is
an integer~$g$ coprime to~$n$ such that $(s_1,\dots,s_\ell)$ is the
sequence of the \lpr/ for the equivalent sequence~$g\alpha$. This is
our main result.

\begin{thm}\label{major}
Each \zfs/ of \length/ in the cyclic group~$\Zn/$ is equivalent to a
sequence whose sum of the least positive representatives is less
than~$n$.
\end{thm}

Such a conclusion does not hold in general for shorter sequences
in~$\Zn/$. Zero-free sequences with lengths at most~$n/2$ and
failing Theorem~\ref{major} are not hard to find.
Consider for example the following  sequences in~$\Zn/$:
\[
\alpha=2^{n/2-1}3\quad\text{for even~$n\ge 6$}\qquad\text{and}\qquad
\beta=2^{(n-5)/2}3^2\quad\text{for odd~$n\ge 9$.}
\]
Here and further on, multiplicities of sequence terms are indicated
by using exponents; for instance $1^32^23$ denotes the sequence
$(1,1,1,2,2,3)$. Both $\alpha$ and~$\beta$ are \zf/, of
lengths~$n/2$ and~$(n{-}1)/2$, respectively. One can check directly
that for each $g$ coprime to~$n$ the sequences $g\alpha$
and~$g\beta$ have sums of their \lpr/ greater than~$n$.

\section{The index of a long minimal zero sequence}
\label{Index}

Chapman, Freeze and Smith defined the {\em index} of a sequence
in~\cite{ChapmanFreezeSmith}. Given a sequence~$\alpha$ in~$\Zn/$,
its index $Index(\alpha)$ is defined as the minimum of $L(g\alpha)$
over all integers~$g$ coprime to~$n$. (Recall that $L(\omega)$
denotes the sum of the \lpr/ of the sequence~$\omega$.) In terms of
the index, Theorem~\ref{major} can be stated as follows.

\begin{thm}\label{majorindex}
Each \zfs/ of \length/ in~$\Zn/$ has index less than~$n$.
\end{thm}

The index of each nonempty zero sequence in~$\Zn/$ is a positive
multiple of~$n$. A zero sequence in~$\Zn/$ is {\em minimal} if each
of its nonempty proper subsequences is \zf/. The question about the
minimal zero sequences with index exactly~$n$ was studied from
different points of view.

For instance, let $\ell(\Zn/)$ be the minimum integer such that
every minimal zero sequence~$\alpha$ in~$\Zn/$ of length at
least~$\ell(\Zn/)$ satisfies ${Index(\alpha)=n}$. Gao \cite{Gao}
proved the estimates $\left\lfloor(n+1)/2\right\rfloor+1\le
\ell(\Zn/)\le n-\left\lfloor(n+1)/3\right\rfloor+1$ for~$n\ge 8$
($\lfloor x\rfloor$ denotes the greatest integer not exceeding~$x$).
Based on Theorem~\ref{major}, here we determine $\ell(\Zn/)$ for
all~$n$.

The proof comes down to the
observation that each minimal zero sequence of length greater
than $n/2+1$ in $\Zn/$ has index~$n$. Indeed, remove one term~$a$
from such a sequence~$\alpha$; this yields a \zfs/~$\alpha'$ of
\length/. By Theorem~\ref{majorindex}, $Index(\alpha')<n$. Since
$\v{ga}\le n$ for any integer~$g$, it follows that $Index(\alpha)\le
Index(\alpha')+n<2n$. So ${Index(\alpha)=n}$, and we obtain
$\ell(\Zn/)\le \lfloor n/2\rfloor+2$ for all~$n$. Now consider the
following  sequences in~$\Zn/$:
\[
\alpha=2^{n/2-1}3({-}1)\ \text{for even~$n\ge
6$}\quad\text{and}\quad \beta=2^{(n-5)/2}3^2({-}1)\ \text{for
odd~$n\ge 9$.}
\]
These modifications of the examples at the end of the previous
section show that the upper bound $\ell(\Zn/)\le \lfloor
n/2\rfloor+2$ is tight for even $n\ge 6$ and odd $n\ge 9$. Indeed,
$\alpha$ and~$\beta$ are minimal zero sequences, of respective
lengths~$n/2+1$ and~$(n+1)/2$. In both cases the length
equals~$\lfloor n/2\rfloor+1$. By the conclusion from the last
paragraph of Section~\ref{TheMainResult}, each of~$\alpha$
and~$\beta$ has index greater than~$n$. (In fact
$Index(\alpha)=Index(\beta)=2n$.)

For the values of $n$ not covered by these examples, that is
$n=2,3,4,5,7$, it is proven in~\cite{ChapmanFreezeSmith} that
$\ell(\Zn/)=1$. It remains to summarize the conclusions.

\begin{prop}\label{index}
The values of $\ell(\Zn/)$ for all~$n>1$ are: If
$n\not\in\{2,3,4,5,7\}$ then $\ell(\Zn/)=\lfloor n/2\rfloor+2$; if
$n\in\{2,3,4,5,7\}$ then $\ell(\Zn/)=1$.
\end{prop}

\section{The maximum multiplicity of a term}
\label{MaxMultTerm}

An extensively used result of Bovey, Erd\H{o}s and Niven
\cite{BoveyErdosNiven} states that each \zfs/ of \lengthell/
in~$\Zn/$ contains a term of multiplicity at least~$2\ell-n+1$. The
authors remark that this estimate is best possible whenever
$\threshold/\le\ell<n$. An improvement for the more interesting
range $n/2<\ell\le \threshold/$ is due to Gao and Geroldinger
\cite{GaoGeroldinger} who showed that $2\ell-n+1$ can be replaced
by~$\max(2\ell-n+1,\ell/2-(n{-}4)/12)$ (for $\ell\ge (n+3)/2$). Here
we obtain a sharp lower bound for each length~$\ell$ greater
than~$n/2$.

Let~$M$ be the maximum multiplicity of a term in a \zfs/ $\alpha$
with \lengthell/ in~$\Zn/$. Clearly $M$ has the same value for all
sequences equivalent to~$\alpha$, and also for the respective
sequences of \lpr/. Therefore, by Theorem~\ref{major}, one may
assume that $\alpha$ is a \pins/ of \lengthell/ and sum $S\le
n{-}1$. Let $\alpha$ contain $u$~ones and $v$~twos. Then
\begin{equation*}
n{-}1\ge S\ge u+2(\ell-u)=2\ell-u,\ n{-}1\ge S\ge
u+2v+3(\ell-u-v)=3\ell-2u-v.
\end{equation*}
These yield $u\ge 2\ell-n+1$ and $2u+v\ge 3\ell-n+1$, respectively.
Since $M\ge\max(u,v)$, it follows that $M\ge
\max\left(2\ell-n+1,\boundnew/\right)$. Now, $2\ell-n+1\ge
\boundnew/$ if and only if $\ell\ge\threshold/$, so two cases arise.

For $\threshold/\le\ell<n$, the lower bound $M\ge 2\ell-n+1$ is best
possible, as already remarked in~\cite{BoveyErdosNiven}. Indeed,
$\alpha=1^{2\ell-n+1}2^{n-\ell-1}$ is a well-defined \pins/ whenever
$n/2<\ell<n$ (note that the last inequality implies $n>2$). It has
length~$\ell$ and sum~$n{-}1$. If in addition $\threshold/\le\ell<n$
then $2\ell-n+1$ is the maximum multiplicity of a term in~$\alpha$,
so $M=2\ell-n+1$.

If $n/2<\ell\le\threshold/$, the lower bound $M\ge \boundnew/$ is
best possible. To show that the equality can be attained, consider
the sequence
\[
\alpha=1^{\ell-\lfloor(n-1)/3\rfloor}2^{\ell-\lfloor(n-1)/3\rfloor}3^{2\lfloor(n-1)/3\rfloor-\ell}.
\]
It is well defined unless $n$ is divisible by~3 and $\ell=2n/3-1$;
this case will be considered separately. The multiplicities of 1, 2
and~3 are nonnegative integers for all other
values of~$n$ and~$\ell$ satisfying
$n/2<\ell\le\threshold/$ (which also
implies $n>3$). So $\alpha$ is a \pins/ with length~$\ell$,
sum~$3\lfloor(n{-}1)/3\rfloor\le n-1$ and two terms of maximum
multiplicity which equals~$\boundnew/$. In the exceptional case
mentioned above, the example $\alpha=1^{n/3}2^{n/3-1}$ shows that
$M=\boundnew/$ is attainable, too.

We proved the following tight piecewise linear lower bound.
\begin{prop}\label{multterm}
Let $n$ and $\ell$ be integers satisfying $n/2<\ell<n$. Each \zfs/
of length~$\ell$ in~$\Zn/$ has a term with multiplicity:
\begin{enumerate}
\item[a)] at least $2\ell-n+1$ if $\threshold/\le\ell<n$;

\item[b)] at least $\boundnew/$ if $n/2<\ell\le\threshold/$.
\end{enumerate}
These estimates are best possible.
\end{prop}

Essentially speaking, the arguments above yield an explicit
description of the \zfs/s in~$\Zn/$ with a given \lengthell/. This
description is included in Section~\ref{ConcludingRemarks}. Here we
only note that the equality $M=\max(u,v)$ holds for each \pins/
$\alpha$ of \length/ and sum at most~$n{-}1$. Indeed, fix
$2\ell-n+1$ ones in~$\alpha$ (this many ones are available in view
of $u\ge 2\ell-n+1$). The remaining part~$\alpha'$ has length
$n-1-\ell$ and sum~$\le 2(n-1-\ell)$, so the average of its terms is
at most~2. It readily follows that $\alpha'$ contains at least as
many ones as terms greater than~2.

\section{The maximum multiplicity of a generator}
\label{MaxMultGenerator}

Given a \zfs/ in~$\Zn/$, what can be said about the number of
generators it contains? As usual, here a {\em generator} means an
element of~$\Zn/$ with order~$n$. This question attracted
considerable attention and effort, for sequences of \length/. Even
the existence of one generator in such a sequence (which follows directly from Theorem~\ref{main}) does not seem
immediate. It was proven by Gao and
Geroldinger~\cite{GaoGeroldinger}. Improving on their result,
Geroldinger and Hamidoune~\cite{GeroldingerHamidoune} obtained the
following theorem. A \zfs/~$\alpha$ of length at least~$(n{+}1)/2$
in~$\Zn/$ ($n\ge 3$) contains a generator with multiplicity~3 if $n$
is even, and with multiplicity $\left\lceil(n{+}5)/6\right\rceil$ if
$n$ is odd ($\lceil x\rceil$ denotes the least integer greater than
or equal to~$x$). These bounds are sharp if $\alpha$ ranges over the
\zfs/s in~$\Zn/$ of {\em all} lengths~$\ell\ge(n{+}1)/2$.

On the other hand, the above estimates do not reflect the length
of~$\alpha$. One can be more specific by finding best possible
bounds for each length~$\ell$ in the range~$\left(n/2,n\right)$.

Denote by~$m$ the maximum multiplicity of a generator in a \zfs/
$\alpha$ with \lengthell/ in~$\Zn/$. By Theorem~\ref{major}, we may
assume again that $\alpha$ is a \pins/ of \lengthell/ and sum at
most~$n{-}1$; the point of interest now is the maximum
multiplicity~$m$ of a term coprime to~$n$. Let $\alpha$ contain
$u$~ones and $v$~twos, as in Section~\ref{MaxMultTerm}. It was shown
there that $u\ge 2\ell-n+1$, and because 1 is coprime to~$n$, we
have $m\ge 2\ell-n+1$.

If $n$ is even, the sequence $1^{2\ell-n+1}2^{n-\ell-1}$ shows that
this bound is sharp.

If $n$ is odd then 2 is coprime to~$n$, so $m\ge\max(u,v)$. But if
$M$ is the maximum multiplicity of a term in~$\alpha$ then $m\le M$,
and also $M=\max(u,v)$ by the remark after
Proposition~\ref{multterm}. Hence $M=m$, so the answer in the case
of an odd~$n$ coincides with the one from the previous section.

The conclusions are stated in the next proposition.

\begin{prop}\label{multgenerator}
Let $n$ and $\ell$ be integers satisfying $n/2<\ell<n$, and let
$\alpha$ be a \zfs/ of length~$\ell$ in~$\Zn/$.
\begin{enumerate}
\item[a)] For $n$  even, $\alpha$ contains a generator of multiplicity at
least $2\ell-n+1$. This estimate is best possible.

\item[b)] For $n$  odd, $\alpha$ contains a generator of multiplicity at
least $2\ell-n+1$ if $\threshold/\le\ell<n$, and at least
$\boundnew/$ if $n/2<\ell\le\threshold/$. These estimates are best
possible.
\end{enumerate}
\end{prop}

The theorem of Geroldinger and Hamidoune~\cite{GeroldingerHamidoune}
can be regarded as an extremal case of
Proposition~\ref{multgenerator}, obtained by setting $\ell=n/2+1$ if
$n$ is even, and $\ell=(n{+}1)/2$ if $n$ is odd.

\section{A function related to \zfs/s}
\label{NewFunction}

For \pin/s $n$ and $k$, where $n\ge k$, let $h(n,k)\ge k$ be the
least integer
such that each sequence in~$\Zn/$ with at
least $k$ distinct terms and length $h(n,k)$ contains a nonempty
zero sum. The function $h(n,k)$ is a natural analogue of a function
introduced by Bialostocki and Lotspeich~\cite{BialostockiLotspeich}
in relation to the renowned theorem of Erd\H{o}s, Ginzburg
and~Ziv~\cite{ErdosGinzburgZiv}.

It is trivial to notice that $h(n,k)=k$ whenever $k$ is greater than
or equal to the {\em Olson's constant} of the group~$\Zn/$. Olson's
constant $Ol(G)$ of an abelian group~$G$ is the least \pin/~$t$ such
that every subset of~$G$ with cardinality~$t$ contains a nonempty
subset whose sum is zero. Erd\H{o}s
\cite{ErdosGraham} conjectured that $Ol(G)\le\sqrt{2|G|}$ for each
abelian group~$G$; here $|G|$ is the order of~$G$. The best known
upper bound for~$Ol(G)$ is due to Hamidoune and Z\'emor
\cite{HamidouneZemor} who proved that
$Ol(G)\le\left\lceil\sqrt{2|G|}+\gamma(|G|)\right\rceil$, where
$\gamma(n)=O\left(n^{1/3}\log n\right)$. On the other hand, the set
$\{1,2,\dots,k\}$ where $k$ is the greatest integer such that
$1+2+\cdots+k<n$, yields the obvious lower bound $Ol(\Zn/)\ge
\left\lfloor\left(\sqrt{8n-7}-1\right)/2\right\rfloor+1$.

As for values of~$k$ less than~$Ol(\Zn/)$, by using
Theorem~\ref{major} one can determine $h(n,k)$ for
all~$k\le\left(\sqrt{4n-3}+1\right)/2$.

\begin{prop}\label{hnk}
Let $n\ge k$ be \pin/s such that $k\le(\sqrt{4n-3}+1)/2$. Then
\begin{equation*}\label{hnkvalues}
h(n,k)=n-\frac12(k^2-k).
\end{equation*}
\end{prop}

\begin{pf}
The claim is true for $k=1$, so let $k>1$. Denote $\ell=n-(k^2-k)/2$
and notice that $2\le k\le\left(\sqrt{4n-3}+1\right)/2$ is
equivalent to $n/2<\ell<n$. We show that each \zfs/ $\alpha$ of
length~$\ell$ in~$\Zn/$ contains fewer than $k$ distinct terms; then
$h(n,k)\le n-(k^2-k)/2$ by the definition of~$h(n,k)$.

By Theorem~\ref{major} one may regard
$\alpha$ as a \pins/ of length~$\ell$ and sum $S\le n{-}1$. An easy
computation shows that $\alpha$ has at least $2\ell-S$ ones. So
$\alpha=1^{2\ell-S}\beta$, where
$\beta$ is a sequence of length~$S-\ell$ and sum~$2(S-\ell)$. Let
there be $m$ distinct terms in~$1^{2\ell-S}\beta$; then $\beta$ has
$m-1$ distinct terms greater than~1. Because $k>1$, we may assume
$m>1$. Choose one occurrence for each of the $m{-}1$ distinct terms
in~$\beta$ and replace these occurrences by $2,3,\dots,m$. Next,
replace each remaining term by~1. The sum
of~$\beta$ does not increase,
so $ 2(S-\ell)\ge
(2+3+\cdots+m)+(S-\ell-m+1)$.
Combined with $S\le n-1$, this leads to $m^2-m-2(n-\ell-1)\le 0$.
Hence
\[ m\le\frac 12\left(\sqrt{8(n-\ell)-7}+1\right)=\frac
12\left(\sqrt{4(k^2-k)-7}+1\right)<k.
\]
Therefore $2\le k\le\left(\sqrt{4n-3}+1\right)/2$ implies $h(n,k)\le
n-(k^2-k)/2$.

Now consider the sequence $\alpha=1^{\ell-k+1}23\dots k$, where
$\ell=n-(k^2-k)/2-1$. Whenever $2\le
k\le\left(\sqrt{4n-3}+1\right)/2$ and $(n,k)\ne (3,2)$, there are
$k$ distinct terms in~$\alpha$ because these conditions
imply $\ell-k+1\ge 1$. Also $\alpha$ has length~${\ell\ge k}$ and is
\zf/ since the sum of its \lpr/ is~$n{-}1$. It follows that
$h(n,k)\ge n-(k^2-k)/2$. The same lower bound holds for $n=3$, $k=2$
by the definition of~$h(n,k)$. Hence $h(n,k)\ge n-(k^2-k)/2$ for
all~$n$ and~$k$ satisfying $2\le k\le\left(\sqrt{4n-3}+1\right)/2$,
which completes the proof. \qed
\end{pf}

The example $\alpha=1^{\ell-k+1}23\dots k$ in the last proof yields
the
lower bound $h(n,k)\ge n-(k^2-k)/2$ for
$k\le \left(\sqrt{8n-7}-1\right)/2$ which is a weaker constraint
than $k\le\left(\sqrt{4n-3}+1\right)/2$ if $n>7$. So the following
query is in order here.

\begin{quest}
Does the equality
\begin{equation*}
h(n,k)=n-\frac12(k^2-k)
\end{equation*}
hold true whenever $k\le \left(\sqrt{8n-7}-1\right)/2$?
\end{quest}

\section{Concluding remarks}
\label{ConcludingRemarks}

Among other consequences, Theorem~\ref{major} yields various
explicit descriptions of the \zfs/s in~$\Zn/$ with a given
\lengthell/. We include one such description mentioned in
Section~\ref{MaxMultTerm}, skipping over the easy justification.

Let $n$ and $\ell$ be integers satisfying $n/2<\ell<n$. An arbitrary
\zfs/~$\alpha$ of length~$\ell$ in~$\Zn/$ has one of the equivalent
forms specified below.
\begin{enumerate}

\looseness=-1
\item[1.] If $\threshold/\le\ell<n$ then
$\alpha\cong1^u\beta$, where $u\ge 2\ell{-}n{+}1$ and $\beta$ is a
sequence of length~$\ell{-}u$ in~$\Zn/$, without ones and satisfying
$L(\beta)\le n{-}1{-}u$.

\item[2.] If $n/2<\ell\le\threshold/$ there are two possibilities:

\begin{enumerate}
\item[a)] $\alpha\cong1^u\beta$, where $u\ge
\ell/2$ and $\beta$ is a sequence of length~$\ell-u$ in~$\Zn/$,
without ones and satisfying $L(\beta)\le n-1-u$.

\item[b)] $\alpha\cong1^u2^v\beta$, where
\[
u\le\frac{\ell}{2},\quad \min(u,v)\ge 2\ell-n+1,\quad \max(u,v)\ge
\bound/,
\]
and $\beta$ is a sequence of length~$\ell-u-v$ in~$\Zn/$, without
ones and twos and satisfying $L(\beta)\le n-1-u-2v$.
\end{enumerate}
\end{enumerate}

A closer look at the description shows that the structure of the
\zfs/s with lengths~$\ell$ satisfying $n/2<\ell\le \threshold/$ is
significantly more involved than the one for~$\ell$ in the
range~$\threshold/\le\ell<n$ considered in~\cite{Gao}.

Yet another application of the main result concerns zero-sum
problems of a different flavor. Let $n$ and~$k$ be integers such
that $n/2<k<n$. By using Theorem~\ref{major}, one can determine the
structure of the sequences in~$\Zn/$ with length $n-1+k$ that do not
contain $n$-term zero subsequences. Such a characterization in turn
has consequences related to variants of the Erd\H{o}s--Ginzburg--Ziv
theorem \cite{ErdosGinzburgZiv} and deserves separate treatment.
Questions of this kind will be considered in a forthcoming article.


\begin{thebibliography}{00}

\bibitem{BialostockiLotspeich}
A.~Bialostocki, M.~Lotspeich, Some developments of the
Erd\H{o}s--Ginzburg--Ziv theorem. I., in: Sets, Graphs and Numbers
(Budapest, 1991), 97--117, Colloq.\ Math.\ Soc.\ J\'anos Bolyai 60,
North-Holland, Amsterdam, 1992.

\bibitem{BoveyErdosNiven}
J.\ D.\ Bovey, P.\ Erd\H{o}s and I.\ Niven, Conditions for a zero
sum modulo~$n$, Canad.\ Math. Bull.\ 18 (1) (1975), 27--29.

\bibitem{ChapmanFreezeSmith}
S.\ T.\ Chapman, M.\ Freeze and W.\ W.\ Smith, Minimal
zero-sequences and the strong Davenport constant, Discrete Math. 203
(1--3) (1999), 271--277.

\bibitem{ErdosGraham}
P.\ Erd\H{o}s, Problems and results on combinatorial number theory,
in: A Survey of Combinatorial Theory, J.\ N.\ Srivastava et al.
(eds.), North-Holland, Amsterdam, 1973, 117--138.

\bibitem{ErdosGinzburgZiv}
P.\ Erd\H{o}s, A.\ Ginzburg and A.\ Ziv, Theorem in the additive
number theory, Bull.\ Res.\ Council Israel 10F (1961), 41--43.



\bibitem{Gao}
W.\ D.\ Gao, Zero sums in finite cyclic groups, Integers 0 (2000),
A12, 7pp. (electronic).


\bibitem{GaoGeroldinger}
W.\ D.\ Gao, A.\ Geroldinger, On the structure of zerofree
sequences, Combinatorica 18 (4) (1998), 519--527.

\bibitem{GeroldingerHamidoune}
A.\ Geroldinger, Y.\ O.\ Hamidoune, Zero-sumfree sequences in cyclic
groups and some arithmetical application, J.\ Th\'eor.\ Nombres
Bordeaux 14 (1) (2002), 221--239.

\bibitem{HamidouneZemor}
Y.\ O.\ Hamidoune, G.\ Z\'emor, On zero-free subset sums, Acta
Arith.\ 78 (2) (1996), 143--152.

\bibitem{SmithFreeze}
W.\ W.\ Smith, M.\ Freeze, Sumsets of zerofree sequences, Arab.\ J.\
Sci.\ Eng.\ Sect.\ C Theme Issues 26 (1) (2001), 97--105.


\end{thebibliography}
\end{document}